\documentclass[12pt]{amsart}
\usepackage{amscd}

\newtheorem{theorem}{Theorem}
\newtheorem{lemma}[theorem]{Lemma}
\newtheorem{corollary}{Corollary}
\newtheorem{proposition}[theorem]{Proposition}

\theoremstyle{definition}

\newtheorem*{example}{Example}

\newtheorem*{notation}{Notation}

\theoremstyle{remark}
\newtheorem*{remark}{Remark}

\usepackage{amsfonts}
\newcommand{\field}[1]{\mathbb{#1}}
\newcommand{\Q}{\field{Q}}
\newcommand{\R}{\field{R}}
\newcommand{\Z}{\field{Z}}

\newcommand{\A}{\field{A}}

\newcommand{\cc}{\mathcal C}
\newcommand{\cT}{\mathcal T}
\newcommand{\cU}{\mathcal U}

\renewcommand{\a}{\alpha}

\newcommand{\D}{\Delta}

\newcommand{\g}{\gamma}
\newcommand{\G}{\Gamma}

\newcommand{\fg}{\mathfrak g}

\newcommand{\fu}{\mathfrak u}
\newcommand{\fl}{\mathfrak l}

\newcommand{\bs}{\backslash}
\newcommand{\ra}{\rightarrow}
\newcommand{\bp}{\begin{proof}}
\newcommand{\ep}{\end{proof}}

\begin{document}

\title[Unipotent    Generators   for    Arithmetic   Groups]{Unipotent
Generators for Arithmetic Groups}

\author{T. N. Venkataramana}

\address{School of Mathematics, Tata Institute of Fundamental 
Research, Homi Bhabha Road, Bombay - 400 005, INDIA.}

\email{venky@math.tifr.res.in}

\subjclass{Primary  Secondary\\
T. N. Venkataramana,
School of Mathematics, Tata Institute of Fundamental 
Research, Homi Bhabha Road, Bombay - 400 005, INDIA.
venky@math.tifr.res.in}

\date{}

\begin{abstract} We sketch  a simplification of proofs  of old results
on  the arithmeticity  of  the group  generated  by opposing  integral
unipotent radicals contained in higher rank arithmetic groups.
\end{abstract}

\maketitle

\section{Introduction}

A well known theorem of Jacques  Tits \cite{Tits1} says that if $n\geq
3, k \geq 1$ are integers, then the group generated by upper and lower
triangular  unipotent matrices  in the  principal congruence  subgroup
$SL_n(k\Z)$ of level $k$, has finite index in $SL_n(\Z)$. This theorem
admits a generalisation which will be described below (for definitions
of the terms involved, see section \ref{preliminaries}). \\

Let  $G \subset  SL_n$ be  a semi-simple  $\Q$-simple algebraic  group
defined  over  $\Q$  and  let  $Q\subset  G$  be  a  proper  parabolic
$\Q$-subgroup with unipotent radical $U^+$.  Let $U^-$ be the opposite
unipotent radical  and for an integer  $k \geq 1$, denote  by $E_Q(k)$
the subgroup generated by $U^+\cap SL_n(k\Z)$ and $U^-\cap SL_n(k\Z)$.
The aim of this  note is to provide a proof  (which is perhaps simpler
and more  uniform than  the existing  ones in  the literature)  of the
following result.

\begin{theorem} \label{elementary}  If $\R-rank (G) \geq  2$, then the
group $E_Q(k)$ is an arithmetic  subgroup of $G(\Q)$, i.e.  has finite
index in $G(\Z)=G\cap SL_n(\Z)$.
\end{theorem}

\begin{remark} Every  semi-simple $\Q$-simple  algebraic group  $G$ is
the group obtained by (Weil)  restriction of scalars, of an absolutely
simple algebraic group  $\mathcal G$ defined over a number  field $K$ :
$G=R_{K/\Q}(\mathcal   G)$.   Theorem   \ref{elementary}  is   due  to
\cite{Tits1}  If  $\mathcal  G$  is a  Chevalley  group  with  $K-rank
(\mathcal  G)\geq  2$.  For  most  of  the classical  groups,  Theorem
\ref{elementary}  was  proved  by Vaserstein  \cite{Vaserstein1},  and
\cite{Raghunathan1}  proved  it for  general  groups  of $\Q-  rank  $
at least two. The remaining cases were proved in \cite{V1}. \\

These references  prove Theorem  \ref{elementary} when  $Q$ is  a {\it
minimal} parabolic  $\Q$- subgroup; however, as  observed in \cite{V3}
and \cite{Oh2}, the general case follows easily from this case.
\end{remark}

\begin{remark}  A generalisation  of Theorem  \ref{elementary} to  the
case  when the  arithmetic group  $G(\Z)$ is  replaced by  any Zariski
dense  discrete   subgroup  of   $G(\R)$  is  proved   in  \cite{Oh1},
\cite{Benoist-Oh} , \cite{Benoist-Oh2} and \cite{Benoist-Miquel}; the
proofs of these results  are of a very different nature  and we do not
consider this situation. In fact,  the proofs in these references {\it
make use of} Theorem \ref{elementary}.
\end{remark}

\begin{remark}  The proof  of Theorem  \ref{elementary} given  here is
uniform;  however, this  is based  on Theorem  \ref{maintheorem} whose
proof  is not  quite uniform  but works  especially well  (see Section
\ref{proof1}) when the group $G(\R)$ is  {\it not} a product of simple
Lie groups of  real rank one.  In  case $G(\R)$ {\it is}  a product of
rank one groups, a more complicated argument is needed and we give the
proof in Section \ref{proof2}.
\end{remark}

We  now describe  another result  from which  Theorem \ref{elementary}
will be  derived. Let $G$  be a $\Q$-simple  group with $\Q  -rank (G)
\geq 1$.  Let $P\subset G$  be a  proper {\it maximal}  parabolic $\Q$
-subgroup. Denote  by $U^+$ (resp.  $U^-$) the unipotent  radical (the
``opposite  unipotent   radical'')  of  $P$   and  let  $P=LU^+$   be  a
Levi-decomposition of  $P$; set $P^-=LU^-$, the  parabolic subgroup of
$G$ ``opposite'' to $P$. Clearly $L$ normalises $U^{\pm}$. \\

Denote by  $M$ the connected component  of the Zariski closure  of the
group $L(\Z)$ of  integer points of $L$. The group  $M$ is non-trivial
if and only if $\R-rank  (G)\geq 2$ (Lemma \ref{Mpositive}). We denote
by $V^{\pm}$ the commutator group $[M,U^{\pm}]$. \\

For  each  $k\geq  1$,  denote   by  $F(k)$  the  group  generated  by
$V^{\pm}(k\Z)$  and  $M(k\Z)$. Denote  by  $Cl(F(k))$  the closure  of
$F(k)$ in  the group $G(\A  _f)$ of finite  adeles $\A _f$  over $\Q$.
Let $\G_k=G(\Q)\cap Cl(F(k))$.  The group  $\G _k$ has finite index in
$G(\Z)$ (Lemma \ref{congruenceclosure}); it is the smallest congruence
subgroup of $G(\Z)$ containing $F(k)$.  We prove:

\begin{theorem} \label{commutatortheorem} If $\R-rank (G)\geq 2$, then
$F(k)$ contains the commutator subgroup $[\G _k,\G _k]$:
\[ [\G _k,\G _k] \subset F(k). \]  
\end{theorem}

We  now  show  that Theorem  \ref{commutatortheorem}  implies  Theorem
\ref{elementary}.   By  the  Margulis  normal  subgroup  theorem,  the
commutator  $[\G _k,  \G  _k]$ has  finite index  in  the higher  rank
arithmetic  group  $\G  _k$  and is  therefore  an  arithmetic  group;
therefore,    by    Theorem   \ref{commutatortheorem},    the    group
\[F(k)=<V^+(k),V^-(k),  M(k)>\] is  an  arithmetic  group.  Since  $M$
normalises $V^+$ and $V^-$ it follows that $F(k)$ normalises the group
$E'_P(k)=<V^+(k),V^-(k)>$ generated by $V^{\pm}(k)$.  Therefore, again
by  the normal  subgroup theorem,  $E'_P(k)$ is  an arithmetic  group.
Since $E'_P(k)$ is contained  in the group $E_P(k)=<U^+(k),U^-(k)>$ it
follows that the group $E_P(k)$ is  an arithmetic group for every {\it
maximal} parabolic $\Q$-subgroup $P$ of $G$. \\

Let  $Q\subset  G$   be  a  parabolic  $\Q$-subgroup   as  in  Theorem
\ref{elementary}.  Fix a  maximal parabolic  $\Q$-subgroup $P$  of $G$
containing   $Q$.   Then    $U^{\pm}_Q   \supset   U^{\pm}_P$.   Hence
$E_Q(k)=<U^+(k),U^-(k)>$          contains          the          group
$E_P(k)=<U^+_P(k),U^-_P(k)>$. By the  preceding paragraph, $E_P(k)$ is
arithmetic   and   hence  so   is   $E_Q(k)$;   this  proves   Theorem
\ref{elementary}. \\

Theorem  \ref{commutatortheorem}  is  deduced  from a  result  on  the
centrality of  the kernel  for a  map between  two completions  of the
group $G(\Q)$.  This  centrality is somewhat analogous to  that of the
centrality of  the congruence  subgroup kernel  (in the  case $\R-rank
(G)\geq 2$), except that the congruence subgroup kernel is a compact (
profinite) group.  In our case, it is not clear, a priori, that $C$ is
even locally  compact (it will  follow after the  fact that $C$  is in
fact  finite).   The  details  of the  construction  of  the  relevant
completion will  be given  in Section \ref{preliminaries}.  We briefly
describe the construction here. \\

Equip the  group $G(\Q)$ with  the topology $\mathcal T$  generated by
the various cosets $\{gF(k): g\in G(\Q),k \geq 1\}$ where $F(k)$ is as
in  Theorem   \ref{commutatortheorem}.   Then  we  prove   in  Section
\ref{preliminaries}  the following  proposition  (note  that even  in
the proposition   the  assumption  of  higher  real
rank is necessary).

\begin{proposition} \label{topologicalgroup} If  $\R -rank (G)\geq 2$,
then  the group  G,  equipped  with the  topology  $\mathcal  T$ is  a
topological group.
\end{proposition}

The  topological group  $(G,{\mathcal T})$  then admits  a (two-sided)
completion $\widehat G$  which can be shown to  map onto $\overline{G}
\subset G(\A_f)$ where $\overline G$ is  the closure of $G(\Q)$ in the
finite adelic  group $G(\A _f)$ (  also referred to as  the congruence
completion of  $G(\Q)$); if  $G$ is simply  connected, then  by strong
approximation, ${\overline G}=G(\A _f)$. Then we get an exact sequence
\[ 1 \ra C  \ra {\widehat G} \ra {\overline G}  \ra 1\] of topological
groups; the map ${\widehat G} \ra {\overline G}$ can be shown to be an
open map. The kernel $C$ is closed in $\widehat G$; however, it is not
clear a-priori that $C$ is even compact). The main result of the paper
is

\begin{theorem} \label{maintheorem} If $\R-rank  (G) \geq 2$, then the
kernel $C$ is central in $\widehat G$.
\end{theorem}

It   can   now  be   seen   why   this  centrality   implies   Theorem
\ref{commutatortheorem}. By  the definition of  the group $\G  _k$ (as
the smallest congruence subgroup containing $F(k)$), the groups $F(k)$
and $\G _k$ have the same closure  in $G(\A _f)$.  The openness of the
map ${\widehat  G} \ra {\overline  G}$ then implies that  the closures
${\widehat \G _k}, {\widehat F(k)}$ in  $\widehat G$ of the groups $\G
_k, F(k)$ have the property that ${\widehat \G _k} \subset C {\widehat
F(k)}$.    Therefore,    by   the    centrality   of    $C$   (Theorem
\ref{maintheorem}) we see that
\[[\G _k,  \G _k] \subset  [{\widehat \G _k},{\widehat \G  _k}] \subset
[{\widehat  F(k)},{\widehat  F(k)}]  \subset   {\widehat  F(k)}  .  \]
Therefore, we get
\[[\G_k,  \G _k]  \subset G(\Q)\cap  {\widehat F(k)}.   \] From  Lemma
\ref{intersection}  we  have  that $G(\Q)\cap  {\widehat  F(k)}=F(k)$;
therefore   we   get   $[\G_k,\G_k]\subset  F(k)$,   proving   Theorem
\ref{commutatortheorem}. \\

The  centrality of  $C$ is  deduced in  Section \ref{proof1}  when the
group $M$ is not abelian; this is  shown to be a simple consequence of
strong approximation.  When the group $M$ is not abelian, the proof is
more complicated  and is is  dealt with in Section  \ref{proof2}; this
involves the  analogues of some  results which are  essentially proved
(but stated  only for the  congruence subgroup kernel, instead  of our
group $C$) in \cite{V2}. \\

{\bf Acknowledgments} I  thank the organisers for inviting  me to take
part in the conference and to contribute an article to the proceedings
of the  conference. The support of  JC Bose fellowship for  the period
2021-2025 is gratefully acknowledged.

\newpage

\section{Preliminaries} \label{preliminaries}

\subsection{Topological Groups} \label{topgroups}

Let  $G$  be a  group  and  $\cc=\{W\}$  a (countable)  collection  of
subgroups $W$ with $\cap _{W\in \cc}  W =\{1\}$, and such that for any
finite set  $F\subset \cc$ there  exists a  subgroup $V \in  \cc$ such
that
\[V \subset  \cap _{W\in  F} W .\]  Let $\cT$ be  the topology  on $G$
generated by the cosets $xW$ with $x\in G$ and $W\in \cc$.

\begin{lemma}  \label{topologicalgroup}  The  pair  $(G,  \cT)$  is  a
topological group if and only if for any $x\in G$ and $W\in \cc$ there
exists a subgroup $V\in \cc$ such that $xWx^{-1}\supset V$.
\end{lemma}

\begin{proof} Suppose for $x\in G$ and $W\in \cc$, there exists $V \in
\cc$ such  that $xWx^{-1}\supset V$. Let  $x,y \in G$ and  put $z=xy$;
let $\cU$ be a neighbourhood of $z$. There exists $W\in \cc$ such that
$zW$ is a neighbourhood of $z$ and $zW\subset \cU$. By our assumptions
on $\cc$,  there exists a  $V \in \cc$  such that $V  \subset yWy^{-1}
\cap W$. We then get
\[  xVyV=xyy^{-1}VyV   \subset  xyWW=xyW=zW,   \]  proving   that  the
multiplication    map   $(x,y)\mapsto    xy=z$   is    continuous   at
$(x,y)$.  Moreover, $x^{-1}W\supset  (Wx)^{-1}=(xx^{-1}Wx)^{-1}\supset
(xV)^{-1}$ for some $V \in \cc$,  proving the continuity of $x \mapsto
x^{-1}$. Therefore, $(G,\cT)$ is a topological group. \\

If the pair  $(G,\cT)$ is a topological group, then  the map $x\mapsto
x^{-1}$ is continuous,  and hence given $W \in \cc$  and $x\in G$, the
group $xWx^{-1}$  is an open subgroup  and hence contains some  $V \in
\cc$.
\end{proof}

\begin{example}  Take  $G\subset  SL_n(\Q)$  to  be  a  $\Q$-algebraic
subgroup  and $\cc$  to  be the  collection,  of principal  congruence
subgroups  $G(k\Z):=G\cap SL_n(k\Z)$  of integral  matrices in  $G\cap
SL_n(\Z)$ congruent to the identity  matrix modulo $k$ with $k\geq 1$.
Then $(G,\cT)$ becomes a topological group and the topology on $G(\Q)$
is the {\it congruence topology}.
\end{example}

\begin{remark}  If the  condition of  Lemma \ref{topologicalgroup}  is
satisfied, we say  that a sequence $x_p$ of elements  of $G$ converges
to an  element $y$, if  given a subgroup  $W\in \cc$, there  exists an
integer  $p(W)$ such  that for  all $p\geq  p(W)$, we  have $x_p^{-1}y
\quad yx_p^{-1} \in W$.
\end{remark}

\subsection{Completions of Topological Groups} \label{completions}

Let  $G$ be  a  topological group  and $\cc=\{W  \}$  a collection  of
subgroups  as  in  \ref{topgroups}.  We   will  say  that  a  sequence
$\{x_n\}_{n\geq 1}$ in $G$ is a (two sided) Cauchy sequence if given a
subgroup $W \in \cc$, there exists an integer $n(W)$ such that
\[x_n^{-1}x_{n+m}\in  W, \quad  x_{n+m}x_n^{-1} \in  W \quad  (\forall
\quad n \geq n(W),  \quad \forall \quad m\geq 1). \]  We will say that
two Cauchy  sequences $\{x_n\},\{y_n\}$ are {\it  equivalent} if given
$W\in \cc$, there exists an integer $n(W)$ such that
\[ x_n^{-1}y_n\in W, \quad y_nx_n^{-1}\in W \quad \forall \quad n \geq
n(W).\]  Denote by  $\widehat G$  the  set of  equivalence classes  of
Cauchy sequences; elements of the original group $G$ may be thought of
as the  set of  constant Cauchy  sequences. The  resulting map  $G \ra
\widehat G$ is an embedding (this follows from the assumption that the
intersection $\cap  W_{W\in \cc}=\{1\}$). If  $x=\{x_n\},y=\{y_n\} \in
\widehat G$  denote by  $xy$ and  $x^{-1}$ respectively  the sequences
$\{x_ny_n\}$  and $\{x_n^{-1}\}$;  it  is routine  to  see that  these
sequences  are Cauchy  and we  then get  the structure  of a  group on
$\widehat G$. \\

Given $W \in \cc$, denote by  $\widehat W$ the set of Cauchy sequences
$x=(\{x_n\}$ such  that for  some integer  $n(W)$ we  have $x_n  \in W
\quad \forall  \quad n\geq n(W)$. Then  $\widehat W$ is a  subgroup of
$\widehat  G$.  Write  ${\widehat  \cc}$ for  the  collection of  sets
$\{{\widehat W}: W \in \cc\}$, and by ${\widehat \cT}$ the topology on
$\widehat G$ generated  by the the collection  of cosets $\{x{\widehat
W}:  x\in {\widehat  G}, \quad  W\in  \cc\}$.  Then  $\widehat G$  and
${\widehat \cc}$  satisfy the conditions of  \ref{topgroups} and hence
$({\widehat G}, {\widehat  \cT})$ is a topological  group, referred to
as the  (two sided)  {\it completion}  of $(G,\cT)$  (see \cite{Bour},
Chapter III, Section  3, Exercise 6). It follows  from the definitions
that ${\widehat G}$  is complete in the sense that  Cauchy sequences in
$({\widehat G},{\widehat \cT})$ converge to an element of $\widehat G$.

We first note an easy consequence of the definitions.  

\begin{lemma} \label{intersection} If  $W, G$ is as before 
then the intersection ${\widehat W}\cap G =W$.
\end{lemma}

\bp Let $g=(g,g,g,\cdots)\in  G$ and suppose it  lies in $\widehat
W$. Therefore,  there exists a Cauchy  sequence $\{x_n\}_{n\geq 1}$
such that $x_n  \in W$ for large enough $n$  which is equivalent to
the  constant sequence  $g$.   Therefore, for  $m$  large enough,  the
elements $gx_m^{-1},  x_m^{-1}g$ lie in  $W'$ for any given  $W' \in \cc$; in
particular,  if   we  take   $W'=W$,  we  then see   that  for   $m$  large,
\[g=(gx_m^{-1})x_m \in W.\]
\ep

\begin{notation} Given elements $x,y$ of a group $\G$, we write  
\[^x(y)=xyx^{-1},\quad [x,y]=xyx^{-1}y^{-1}.\] If $\D \subset \G$ is a
subgroup,  we write  $^x(\D)=x \D  x^{-1}$.  If $A,B  \subset \G$  are
subgroups,  then  $[A,B]$  denotes   the  subgroup  generated  by  the
commutators $[a,b]$ with $a\in A \quad b \in B$.
\end{notation}

\begin{example}  Let  $G\subset  SL_n$  be a  linear  algebraic  group
defined over $\Q$. Consider the collection $G(k\Z)=G\cap SL_n(k\Z)$ of
congruence  subgroups  in  $G(\Q)$.   This  collection  satisfies  the
hypotheses of  \ref{topgroups} and hence  we get a  completion denoted
${\overline  G}$  of  $G(\Q)$,  referred to  as  the  {\it  congruence
completion} of $G(\Q)$.  This is also the closure  of $G(\Q)$ embedded
as  a subgroup  of $G(\A  _f)$ where  $\A _f$  is the  ring of  finite
adeles. \\

We  may  also   consider  the  collection  of   subgroups  of  $G(\Q)$
commensurable to  $G(\Z)$ (these  are refereed  to as  {\it arithmetic
subgroups} of  $G(\Q)$); the  collection of arithmetic  subgroups also
satisfy  the  hypotheses  of  \ref{topgroups}.   We  therefore  get  a
completion  ${\widehat  G}_a$ of  $G(\Q)$,  referred  to as  the  {\it
arithmetic completion} of $G(\Q)$. We  then have a surjective open map
${\widehat G}_a  \ra {\overline G}$  of topological groups  split over
$G(\Q)$; the kernel  $C_G$ is seen to be a  profinite (compact) group,
called {\it the congruence subgroup kernel}. \\

The foregoing facts are well known (\cite{BMS}). 
\end{example}

\subsection{Isotropic Algebraic Groups over  $\Q$} \label{isotropic}

In what  follows, $G \subset  SL_n$ is a $\Q$-simple  linear algebraic
group defined  over $\Q$.  It  is said  to be $\Q$-isotropic  if there
exists a  torus $\Q$-isomorphic to the  multiplicative group ${\mathbb
G}_m$ embedded  in $G$; let $S$  in $G$ be a  maximal $\Q$-split torus
and $\Phi=  \Phi (  \fg,S)$ be  the roots  (characters of  $S$ written
additively) of $S$ occurring in the Lie algebra $\fg$ of $G$ under the
adjoint action of $S$. Denote the {\it root space} of $\a \in \Phi$ by
$\fg _\a$, the  subspace of $\fg$ on  which the torus $S$  acts by the
character  $\a \in  \Phi$. Fix  a positive  system of  roots $\Phi  ^+
\subset \Phi$; then $\Phi =\Phi ^+  \cup (-\Phi ^+)$. We write $\a >0$
if $\a \in \Phi ^+$. \\

Denote by $P_0$ the connected subgroup of $G$ whose Lie algebra is the
direct sum $\fg  ^S \oplus _{\a \in  \Phi ^+} \fg _\a$  where $\fg ^S$
denotes the subspace of vectors in  $\fg$ fixed by the split torus $S$
(the Lie algebra  of the centraliser of  $S$ in $G$). Then  $P_0$ is a
minimal parabolic $\Q$-subgroup of $G$.  Let $P \subset P_0$ be a {\it
maximal}  parabolic  $\Q$-subgroup  of  $G$ and  $U^+$  its  unipotent
radical with Lie  algebra $\fu ^+$. Then $\fu ^+  \subset \oplus _{\a>0}
\fg _\a$  and is a  sum $\fu  ^+=\oplus _{\a \in  X} \fg _\a$  of root
spaces for some subset $X\subset \Phi  ^+$ of positive roots. There is
a decomposition (the Levi decomposition) of $P$ as a product $P=LU^+$,
where $L$  is a connected  subgroup of  $G$ containing $S$,  whose Lie
algebra is the direct sum of the root spaces $\fg _{\pm \a}$ with $ \a
\notin X$ and $\fg ^S$. \\

Let $\fu ^-= \oplus _{\a \in X}  \fg _{-\a}$ and $U^-$ the connected (
in fact  unipotent) subgroup of $G$  with Lie algebra $U^-$.   This is
called the  {\it opposite} of $U$;  the group $P^-=U^-L$ is  a maximal
parabolic  $\Q$-subgroup  called  the  {\it  opposite  of  $P$}.   The
multiplication map  $U^-\times P  \ra G$ given  by $(v,p)  \mapsto vp$
identifies the product space $U^-\times P$ as a Zariski dense open set
${\mathcal U}=U^-P$ in $G$ defined over  $\Q$. If $H\subset SL_n$ is a
$\Q$-subgroup, and $\Z_p$  is the ring of $p$-adic  integers, we write
$H(k\Z_p)$ for the subgroup of  elements of $H \cap SL_n(\Z_p)$ viewed
as $n\times  n$-matrices which  are congruent  to the  identity matrix
modulo $k$.
 
\begin{lemma}Let  $k \geq  1$  be  an integer  and  for  a prime  $p$,
consider  the  set  ${\mathcal  U}(k\Z_p)=U^-(k\Z_p)P(k\Z  _p)$  where
$\Z_p$ is the  ring of $p$-adic integers. There exists  a compact open
subgroup $K_p(k)$ of $G(\Z_p)$ contained in ${\mathcal U}(k\Z_p) $.
\end{lemma}

\begin{proof}  The set  ${\mathcal  U}(k\Z_p)$ is  an  open subset  of
$G(\Z_p)$ containing  $1$. A  fundamental system of  neighbourhoods of
identity in $G(\Z_p)$  is given by open subgroups and  hence the lemma
follows.
\end{proof}

\subsection{The Groups $M$, $V^+$ and $V^-$} 

The  group $L(\Z)$  of  integer points  of the  Levi  subgroup $L$  of
subsection  \ref{isotropic} is  not Zariski  dense in  $L$ (since  the
$\Q$-split central  torus of $L$ has  only a finite number  of integer
points). Denote by $M$ the  connected component of the Zariski closure
of $L(\Z)$. The group  $M$ does not change if we  replace $L(\Z)$ by a
finite index subgroup. Since $L(\Q)$ commensurates $M(\Z)$, it follows
that its Zariski closure $L$ normalises $M$.

\begin{lemma} \label{Mpositive}  The dimension  of $M$ is  positive if
and only if $\R-rank (G)\geq 2$.
\end{lemma}

\begin{proof}  By  definition,  $M$  is  the  connected  component  of
identity of  the Zariski  closure of $L(\Z)$;  hence its  dimension is
zero if and only if $L(\Z)$ is  finite. Since $L$ is the Levi subgroup
of  a maximal  parabolic  $\Q$-  subgroup $P$  of  $G$,  we may  write
$L=S_1L_1L_2$ as a product where $S_1$ is a one dimensional $\Q$-split
torus, $L_1$ is  the product of $\Q$-isotropic simple  factors of $L$,
and $L_2$  is a  $\Q$-anisotropic group.  (To see  this, we  write the
semisimple part $L^{ss}$ of $L$ as a product $L_1L_3$ where $L_1$ is a
product of $\Q$-isotropic  simple groups $L_1$ and $L_3$  is a product
of $\Q$-anisotropic  simple groups; then  $L=S_1T_1L^{ss}$ where
$S_1$  is $\Q$  split torus  in  the centre  of  $L$, and  $T_1$ is  a
$\Q$-anisotropic part of the centre of $L$. We may take $L_2=L_3T_1$). \\

If $L(\Z)$  is finite,  then $L_2(\Z)$  is finite  and since  $L_2$ is
$\Q$-anisotropic,  by  the  Godement criterion,  $L_2(\R)/L_2(\Z)$  is
compact and hence  $L_2(\R)$ is compact and has real  rank zero. Since
$L_1(\Z)$ is also finite but $L_1$ is a product of $\Q$ simple groups,
it follows that $L_1$ is trivial  and hence $L=S_1L_2$ where $L_2$ has
real rank $0$.  Consequently the real rank of $L$  is the dimension of
$S_1$ which is one and hence $\R-rank (G)=\R-rank (L)=1$. \\

Conversely,  if  the  real  rank  of   $G$  is  one,  then  the  group
$L(\R)=S_1(\R)L_1(\R)L_2(\R)$  has real  rank one  and hence  $L_1$ is
trivial  and  $L_2$  is  anisotropic over  $\R$;  hence  $L_2(\R)$  is
compact. Therefore, $L(\Z)\simeq  S_1(\Z)L_1(\Z)=\{\pm 1\} L_2(\Z)$ is
finite and hence $M$ has dimension zero.
\end{proof}

Denote  by   $V^{\pm}$  the  group  $[M,U^{\pm}]$   generated  by  the
commutators $mum^{-1}u^{-1}$ with $m\in M$  and $u \in U^{\pm}$. Since
$M$ is normal in $L$ (and  $U^{\pm}$ are normalised by $L$, it follows
that $V^{\pm}$  is normalised by $L$.  It is clear that  $V^{\pm}$ are
unipotent subgroups of $U^{\pm}$.

\begin{lemma} \label{G'} Suppose $G$ is $\Q$-simple and $\Q$-isotropic.\\

[1] The group $U^{\pm}$ normalises $V^{\pm}$. \\

[2] If $\R-rank (G)\geq  2$ then $G$ is generated
by the groups $V^+,V^-$ and $M$.
\end{lemma}

\begin{proof} The action of the reductive group $M$ on the Lie algebra
$\fu ^{\pm}$  is completely  reducible. Consequently, the  lie algebra
$\fu$ splits into the space $(\fu ^{\pm})^M$ of $M$ invariants and the
space of {\it  non-invariants} i.e.  the span of  $mXm^{-1}-X$ with $X
\in \fu ^{\pm}$ and $m \in M$. Since the non-invariants all lie in the
Lie algebra $Lie (V ^{\pm})$, it  follows that $\fu ^{\pm} =(\fu ^{\pm
})^M  + (Lie  V)^{\pm}$  (it is  possible that  the  Lie algebra  {\it
generated} by the non-invariants picks up invariant vectors; hence the
sum may not be  direct).  If $X \in (\fu ^{\pm})^M$, $m  \in M$ and $Y
\in \fu ^{\pm}$, then we have
\[  [X,m(Y)-Y]= [m(X),m(Y)]-[X,Y]=m([X,Y])-[X,Y]  \in (Lie  V)^{\pm}.\]
Therefore,  $U^{\pm}$ normalises  $(Lie V  )^{\pm}$ proving  the first
part. \\

Let $G'$ be the group generated  by $V^{\pm}$ and $M$; since all these
groups are connected, so is $G'$ ; let $\fg '$ be its Lie algebra.  We
will show  that $\fg$  normalises $\fg'$;  the $\Q$-simplicity  of $G$
then implies  that $\fg'=\fg$ and  hence that $G'=G$. Since  the group
$L$ normalises $U^{\pm}$ and also  $M$, clearly $L$ normalises $\fg '$
and hence its Lie algebra $\fl$ normalises $\fg '$.  We therefore need
to check  that $\fu ^{\pm}$  normalises $\fg  '$.  We first  note that
since $M$  is (connected and) normal  in $M$, we have  $m(Z)-Z \in Lie
(M)\subset \fg  '$ if  $Z \in  \fl$.  Since  $[M,U^{\pm}]=V^{\pm}$, it
follows that if $Z \in \fu ^{\pm}$ then $m(Z)-Z \in Lie V^{\pm}\subset
\fg '$. Since  the whole Lie algebra $\fg$ is  spanned by $\fu ^{\pm}$
and $\fl$, we have
\begin{equation}  \label{Mconjugation} m(Z)-Z\in  \fg '  \quad \forall
\quad Z\in \fg. \end{equation}

We have proved in the proof of  the first part of the lemma, that $\fu
^{\pm}=(\fu  ^{\pm  })^M  +  Lie (V^{\pm})$  The  latter  spaces  $Lie
(V^{\pm})$ are  already contained in  $\fg'$.  Therefore, in  order to
verify that the sub-algebras $\fu^{\pm}$ normalise $\fg'$, it is enough
to check that the $M$-invariants in $\fu ^{\pm}$ normalise $\fg'$. \\

Suppose $X\in (\fu  ^+)^M$ and $Y \in  \fu ^{\pm}$.  Fix $m  \in M$. We
compute the bracket
\[[X,m(Y)-Y]=  [X,m(Y)]-[X,Y]=[m(X),m(Y)]-[X,Y]  ,\]  where  the  last
equality  follows because  $X$  is invariant  under  $m\in M$.   Hence
$[X,m(Y)-Y]$    is    $m(Z)-Z$    with   $Z=[X,Y]$.     By    equation
(\ref{Mconjugation}), the bracket $[X,m(Y)-Y] =m (Z)-Z \in \fg '$.  We
have thus  proved that $[ (\fu  ^+)^M, Lie (V^{\pm})] \subset  \fg '$.
If $Z\in Lie (M)$, then $[(\fu ^+)^M, Z]=0 \subset \fg'$. Since $\fg'$
is  generated by  $Lie  (V^{\pm})$ and  $Lie (M)$  and  each of  these
spaces,  upon taking  brackets with  elements of  $(\fu ^+)^M$  lie in
$\fg'$  , it  follows  that  $[(\fu ^+)^M,  \fg']  \subset \fg'$:  the
$M$-invariants   in  $\fu   ^+$  normalise   $\fg'$.   Similarly   the
$M$-invariants in $\fu ^-$ also normalise $\fg '$.  By the last remark
of the preceding paragraph, $\fu ^{\pm}$ normalises $\fg'$. \\

On the  other hand  $\fl$ is  contained in  the normaliser  of $\fg'$.
Therefore all of $\fg$ normalises $\fg'$ and the lemma follows.
\end{proof}

If $H \subset  SL_n$ is an algebraic $\Q$-subgroup,  we write $H(k\Z)$
for the intersection $H\cap SL_n(k\Z)$, where $SL_n(k\Z)$ is the group
of $n\times n$ matrices in $SL_n(\Z)$ congruent to the identity matrix
modulo $k$. Denote by $F(k)$ the group generated by $V^{\pm}(k\Z)$ and
$M(k\Z)$. We note a corollary of lemma \ref{G'}.

\begin{corollary}  \label{F(k)dense} If  $\R -rank  (G) \geq  2$, then
$F(k)$ is Zariski dense in $G$.
\end{corollary}

\begin{proof} Since  $V^{\pm}$ are unipotent $\Q$-groups,  it is clear
that  $V^{\pm}(k\Z)$   are  Zariski  dense  in   $V^{\pm}$.  Moreover,
$M(k\Z)\simeq L(k\Z)$ is Zariski dense  in $M$. Therefore, the Zariski
closure of $F(k)$ contains $V^{\pm}$  and $M$.  By Lemma \ref{G'}, The
Zariski closure of $F(k)$ is equal to $G$.
\end{proof}

\subsection{Strong Approximation}

We recall some well known  results on strong approximation. Suppose $H
\subset SL_n$ be a simply connected semi-simple $\Q$ -simple algebraic
group with $\R-rank  (G) \geq 1$. We work with  the fixed embedding $H
\subset SL_n$. Set $H(\Z)=H  \cap SL_n(\Z)$. Then strong approximation
says  that $H(\Z)$  is  dense  in the  group  $H({\widehat Z})$  where
${\widehat \Z}=\prod  \Z _p$  where $p$ runs  through all  primes. Let
$a,b$  be coprime  integers  and $H(a\Z)=H(\Z)\cap  SL_n(a \Z)$  where
$SL_n(a  \Z)$ are  integral matrices  in $SL_n(\Z)$  congruent to  the
identity matrix modulo the integer $a$.  Then, as before, $H(a \Z)$ is
called the principal congruence subgroup of level $a$.

\begin{lemma}  \label{strongapproximation}   If  $a,b$   are  co-prime
integers, and $H\subset  SL_n$ is a simply  connected $\Q$-simple with
$\R -rank (H)\geq 1$, then $H(a\Z)$ and $H(b\Z)$ generate $H(\Z)$. \\

The same conclusion holds if  $H$ is semisimple, simply connected $\Q$
group which is {\it product}  of $\Q$-simple (simply connected) groups
$H_i$ with $\R-rank (H_i)\geq 1$ for each $i$.

\end{lemma}

\bp This  lemma and the  proof are  well known consequences  of strong
approximation.  For the  sake of  completeness of  the exposition,  we
recall the proof. \\

The  definitions imply  that $H(a\Z)=H(\Z)\cap  H(a{\widehat \Z})$  is
dense in  $H(a {\widehat \Z})$.  Thus  the group $\G ^*$  generated by
the  two  principal  congruence  groups $H(a  \Z),H(b\Z)$  is  a  also
congruence group dense  in the group $H^*$  generated by $H(a{\widehat
\Z})=\prod _pH(a  \Z _p)$  and $H(b{\widehat \Z})=\prod  _p H(b\Z_p)$;
since $a,b$ are coprime, {\it at  each prime $p$}, the group generated
by $H(a\Z_p)$ and  $H(b\Z_p)$ is $H(\Z _p)$ and hence  $H^*$ is all of
$H({\widehat \Z})$.  If  two congruence subgroups of  $H(\Z)$ have the
same  closure in  the congruence  completion $H({\widehat  \Z})$, then
they are  the same; hence  $\G ^*=H(\Z)$ and  first part of  the lemma
follows. \\

The second  part readily follows from  the first part applied  to each
$H_i$.
\ep

\begin{lemma}  \label{congruenceclosure}  The   intersection  $\G  _k=
G(\Q)\cap {\overline F(k)}$  where $\overline F(k)$ is  the closure of
the group $F(k)$ in the congruence completion $\overline G$ of $G(\Q)$
is  an  arithmetic  group  (called the  {\it  congruence  closure}  of
$F(k)$).
\end{lemma}

\bp 

A theorem  of Nori  and Weisfeiler (\cite{Nori},  \cite{W}) says  , in
particular, that  if $\G  \subset G(\Z)$ is  a Zariski  dense subgroup
(and $G$ is $\Q$-simple, $\Q$-isotropic),  then the closure of $\G$ in
the congruence  completion (in this  case $G(\A_f)$) is open).   It is
not difficult to  extend this to the case when  $G$ is not necessarily
simply connected (but  the congruence completion $\overline  G$ of $G$
may not  be all of $G(\A_f)$).  Thus the intersection of  $G(\Q)$ with
the closure of $\G$ is  a congruence (arithmetic) subgroup of $G(\Q)$;
it is the smallest congruence  subgroup of $G(\Z)$ containing $\G$ and
is called the congruence closure of $\G$. \\

Applying  this  to  the   Zariski  dense  subgroup  $F(k)$  (Corollary
\ref{F(k)dense}), we see that the  congruence closure $\G_k$ of $F(k)$
has  finite index  in $G(\Z)$  (in fact,  in this  case,one can  prove
directly by a somewhat lengthy argument  that the closure of $F(k)$ in
$\overline G$ is open, without using Nori-Weisfeiler).

\ep

\begin{remark}

Consider the  group $P_F^{\pm}=MV^{\pm}$  where $V^{\pm}=[M,U^{\pm}]$.
We have seen (Lemma \ref{G'})  that $U^{\pm}$ normalises $V^{\pm}$; it
then follows that $U^{\pm}$ normalises $MV^{\pm}$ as well, since
\[^u(mv)=[u,m]m ^u(v)=m[^{m^{-1}}(u),m] ^u (v)\in MV.\] The groups $M$
and  $V^{\pm}$ are  normalised also  by $L$;  hence $P^{\pm}=LU^{\pm}$
normalises $P_F^{\pm}$. . Since $P_F^{\pm}$ is the semi-direct product
of the $\Q$ groups $M$ and $V^{\pm}$, it follows from definitions that
for varying integers $k$, $M(k)V^{\pm}(k)$  is a fundamental system of
congruence subgroups  of $P_F(\Q)$;  in particular given  $p\in P(\Q)$
and  an integer  $k \geq  1$, there  exists an  integer $l$  such that
$p(M(k)V^{\pm}(k))p^{-1}\supset M(l)V^{\pm}(l)$.

\end{remark}

\begin{lemma}  \label{finitedense} Given  a  Zariski  dense subset  $D
\subset U^-(\Q) $ and an integer $k \geq 1$, there exists a finite set
$F$ in $D$  and an integer $l\geq 1$ such  that the group $M(l)V^-(l)$
is contained in $B$ where $B$ is the group generated by the conjugates
$^v(M(k)):=v(M(k))v^{-1}$ as  $v$ runs through elements  of the finite
set $F$.
\end{lemma}

\bp Fix  an element  $v^* \in  D$. Consider  the algebraic  group $V'$
generated by the elements of the form
\[\phi                (u)=^{v^*}([m^{-1},u])=               [v^*m^{-1}
(v^*)^{-1}]v^*umu^{-1}(v^*)^{-1}=^{v^*}\big{(}(m^{-1})^u(m)\big{)}\]
with $u=(v^*)^{-1}v$ varying through the dense set $D'=(v^*)^{-1}D$ as
$v$ varies in $D$ and $m$ varies  in $M$.  Since $D'$ is Zariski dense
in  $U^-$, it  follows  that  this group  $V'$  is  the group  $^{v^*}
([M,U^{-}]=^{v^*}(V^-)=V^-$ since $U^-$ normalises $V^-$. \\

For reasons of dimension, there exists  a finite set of these elements
$v$  such  that the  elements  $\phi  (u)$  generate a  Zariski  dense
subgroup of the unipotent group $U^-$  as $m$ varies in $M(l)$ and $v$
varies  in $F$.   Since  these  finite set  of  elements  $u$ are  all
rational,  by choosing  the  congruence level  $l'$  suitably, we  may
assume that  $\phi (u)$ are  all elements in  $U^-(k)$ for all  $m \in
M(l')$ and  all $v \in  F$. But a  Zariski dense subgroup  of integral
elements in a  unipotent group (namely $V^-$)  contains $V^-(l'')$ for
some    congruence    level     $l''$.     Moreover,    since    $\phi
(u)={v^*}(m)^v(m^{-1})$,  the   group  $^{v^*}(M(k))$   together  with
$V^-(l'')$ generates a congruence subgroup containing $M(l)V^-(l)$ for
some $l$.  \ep

\begin{proposition} \label{conditionF} Assume $G$ is a $\Q$-simple $\Q$
isotropic algebraic  group with $\R  - rank  (G) \geq 2$.  Given $x\in
G(\Q)$ and $k\geq 1$, there exists an integer $l=l(k,x)$ such that
\[ ^x(F(k))\supset F(l).\]
\end{proposition}

\bp  For   every  $\theta   \in  F(k)$  we   have  $^x(F(k))=^{x\theta
}(F(k))$. Since  $F(k)$ is  Zariski dense  in $G$,  we may  assume, by
replacing $x$ by  $x \theta$ if necessary, that  $x \in U^-P={\mathcal
U}$.  Write $x=vp$ accordingly with $v  \in U^-(\Q)$ and $p \in P(\Q)$
with $v=v(x)$ depending algebraically on $x \in {\mathcal U}$. Then,
\[  ^x(F(k))  \supset  ^x  (M(k)V^+(k))  \supset  ^x  (M(k)V^+(k))\cap
M(k)V^-(k) =\]
\[\supset       ^v(M(l_1)V^+(l_1))\cap       M(k)V^-(k)\supset       ^
v(M(l_1)V^+(l_1)\cap  M(l_1)V  ^-(l_1))  ,\] for  some  integer  $l_1$
(since $v\in  P^-(\Q)$ normalises $MV  ^-$). Since $MV\cap  MV^-=M$ we
get:
\[  ^x(F(k))  \supset ^v(M(l_2))\]  for  some  integer $l_2$  with  $v
=v(x)$. Replacing $x$ by  any $x \g $ with $\g \in  F(k)$, we see that
$^x (F(k))\supset  ^{v(x\g )}(M(l_\g))$ for  some integer $l  _\g$. By
Lemma \ref{finitedense}, for some finite  set $F$ of these $\g$'s, the
group generated by the conjugates
\[^v(M(l_2)),^{v(x\g)}(M(l_\g)  \quad  (\g  \in   F),  \]  contains  a
congruence subgroup  of the form  $M(l\Z)V^-(l)$ for some  integer $l$
and therefore,  $^x(F(k))\supset M(l)V^-(l)$ for some  $l$; similarly,
$^x(F(k))\supset  M(l)V^+(l)$  for  some   $l$  and  hence  $^x(F(k))$
contains the group $F(l)$ generated by $V^{\pm}(l)$ and $M(l)$.

\ep

\subsection{The Group $C$}

By (\ref{topgroups})  and by  Proposition \ref{conditionF} we  get the
following. If $G$ is a $\Q$-isotropic $\Q$-simple algebraic group with
$\R  -rank (G)  \geq  2$,  denote by  $\cT$  the  topology on  $G(\Q)$
generated  by   the  cosets  $xF(k):   x  \in  G(\Q),k\geq   1$.  Then
$(G(\Q),\cT)$  gets   the  structure   of  a  topological   group.  By
(\ref{completions}) The topological group  $(G,\cT)$ admits a two sided
completion $({\widehat G}, {\widehat \cT})$. If , as before, $\overline
G  \subset  G(\A  _f)$  denotes the  {\it  congruence  completion}  of
$G(\Q)$, we get a surjective homomorphism ${\widehat G} \ra {\overline
G}$. This proves Proposition \ref{topologicalgroup}. \\

Since  the group  $F(k)$  lies in  $G(k\Z)$  the principal  congruence
subgroup  of $G(\Z)=G  \cap  SL_n(\Z)$  of level  $k$,  and since  the
$G(k\Z):  k\geq 1$  form  a fundamental  system  of neighbourhoods  of
identity, it follows  that any congruence subgroup  of $G(\Q)$ contains
$G(k\Z)$ for some  $k$ and hence contains $F(k)$. Since  the group $\G
_k$ is the smallest congruence subgroup of $G(\Q)$ containing $F(k)$,
it  follows   that  the   $\G  _k$  form   a  fundamental   system  of
neighbourhoods   of   identity   in   $G(\Q)$   for   the   congruence
topology. Since $F(k)$ is dense n $\G _k$, it follows that if $l$ is a
multiple of $k$, then the quotient  set $F(k)/F(l)$ maps onto the {\it
finite} congruence quotient set $\G_k/\G_l$. Taking inverse limits, it
follows that  ${\widehat F(k)}$  maps onto $Cl(\G_k)=Cl(F(k))$  and the
latter is an open subgroup of $\overline G$. Thus the map ${\widehat G}
\ra {\overline G}$ is an {\it open map}, with kernel $C$, say. \\

The  kernel $C$  is the  inverse image  of the  completions ${\widehat
\G_k}$ as $k$ varies. Moreover, the inverse limit of ${\widehat F(k)}$
is trivial. Hence we get
\begin{equation}  \label{Casinverselimit}   C  =\varprojlim  {\widehat
F(k)} \bs  {\widehat \G _k}/{\widehat F(k)}  =\varprojlim F(k)\bs \G_k
/F(k)  .  \end{equation}  Note   that  the  group  $M(\Z)$  normalises
$V^{\pm}(k\Z)$ and  $M(k\Z)$ and hence normalises  $F(k), \G_k$. Since
$C$  is normal  in $\widehat  G$,  it is  normalised by  $G(\Q)\supset
M(\Z)$;  the above  expression (\ref{Casinverselimit})  of $C$  as the
inverse limit of  the double cosets $F(k)\bs  \G_k/F(k)$ respects this
$M(\Z)$ action.

\newpage
\section{When $M$ is not abelian} \label{proof1}

We   now   prove  the   centrality   of   the  kernel   $C$   (Theorem
\ref{maintheorem}) in the case when $M$  is not abelian.  Since $M$ is
connected reductive  and is (the  connected component of  identity of)
the Zariski  closure of  $L(\Z)$, it follows  that $M(\Z)$  is Zariski
dense  in   $M$;  hence  the   commutator  subgroup  $S=[M,M]$   is  a
(non-trivial) semi-simple $\Q$-group with  $S(\Z)$ being Zariski dense
in  $S$.  Let  $S^*$ denote  the simply  connected cover  of $S$.  Let
$S^*=\prod  S_i$ be  a  product of  $\Q$-simple  groups $S_i$.   Since
$S^*(\Z)$ is Zariski dense in $S^*$, we have that $S_i(\Z)$ is Zariski
dense in  each $S_i$,  and hence each  $S_i(\R)$ is  non-compact; i.e.
$\R-rank (S_i)\geq 1$ for each $i$. \\

Consider an element $x$ in  the double coset $C_k=F(k)\bs \G_k /F(k)$.
Since $F(k)$ is  Zariski dense in $G$, we may  choose a representative
$x   \in  \G_k$   with  $x=vp$   with  $v   \in  U^-(\Q)$   and  $p\in
P(\Q)$.  Moreover,  since the  closure  of  $F(k)$ in  the  congruence
topology on $G(\Q)$ is open, we may  choose $x$ so that for all primes
$p$  dividing the  level $k$,  $x$ lies  in the  open neighbourhood  $
U^-(k\Z  _p)P(k\Z_p)$  of identity  in  $G(\Z_p)$.  Thus the  rational
matrices $v,p$  have a common  denominator, say $a$; but  the elements
$v,p$ are integral at all primes $p$ dividing $k$; in other words, $a$
is coprime to $k$. \\

Fix an element  $m \in M(a^N\Z)$ for some large  $N$.  Since the group
$U^-$ is normalised by $M$ and $v\in U^-(\Q)$ has denominator dividing
$a$,  the   commutator  $[m,v]=mvm^{-1}v^{-1}$  is  integral   and  is
divisible  by $k$  at all  primes  $p$ dividing  $k$. Moreover,  since
$(m,v) \mapsto mvm^{-1}v^{-1}$  is a polynomial in the  entries of $v$
and   $m$   with  integer   coefficients,   for   $N$  large   enough,
$mvm^{-1}v^{-1}$  is integral  at all  primes dividing  $a$; in  other
words, $[m,v]\in  V^-(k\Z)$, where, we  recall, $V^{\pm}=[M,U^{\pm}]$.
Similarly, the commutator $[p^{-1},m] \in (MV)(k\Z)$. \\

We now consider the conjugate  $mxm^{-1}$. We have written $x=vp$ with
$v\in U^-(\Q) $ and $p\in P(\Q)$. Hence
\[mxm^{-1}=mvm^{-1}mpm^{-1}=[m,v]vp[p^{-1},m]=[m,v]x[p^{-1},m].     \]
Thus, if  $m\in M(a^N\Z)$, then  from the discussion in  the preceding
paragraph, as an element of  the double coset $C_k=F(k)\bs \G_k/F(k)$,
$mxm^{-1}\in F(k)xF(k)$,  i.e. $mxm^{-1}=x$ as double  cosets In other
words, the group $M(a^N\Z)$ acts  trivially, under conjugation, on the
element $x \in C_k$. \\

The double  coset $x \in C_k$  may be replaced (since  $F(k)$ has open
closure in the congruence completion of $G(\Q)$), by an element $y \in
\G _k$ such that  for each prime $p$ dividing $a$,  the element $y \in
U^-(\Z_p)P(\Z_p)$. In other words, if  $y=v'p'$ is written as a product
of  $v '\in  U^-(\Q),  p'\in  P(\Q)$, then  $v  '\in U^-(\Z_p),  p'\in
P(\Z_p)$.  In  other words, the  elements $v',p'$ are integral  at all
primes dividing  $a$. Therefore, the  common denominator (say  $b$) of
the rational matrices $v',p'$ is  co-prime to $a$. From the conclusion
of the preceding paragraph, the group $M(b^N\Z)$ acts trivially on the
coset representative $y$ of $x$. \\

Since $S^*(\Q)$ acts on $C_k$ via its image $S(\Q)$ in $M(\Q)$, we see
from the last two paragraphs that, both $S^*(a^N\Z)$ and $S^*(b^N(\Z)$
act trivially on the double coset $x \in C_k$, hence so does the group
generated  by these.   By Lemma  \ref{strongapproximation}, the  group
generated by  these subgroups is  $S^*(\Z)$, and hence  $S^*(\Z)$ acts
trivially on each  coset $x$ in $C_k$. Hence  $S^*(\Z)$ acts trivially
on $C_K$ and by taking inverse  limits, it follows that $S^*(\Z)$ acts
trivially on  the kernel $C$.   But all of  $G(\Q)$ {\it acts}  on the
kernel $C$, and the infinite  group $S^*(\Z)$ acts trivially.  In view
of  the simplicity  of  $G(\Q)$  modulo its  centre,  it follows  that
$G(\Q)$,  and hence  ${\widehat  G}$,  act trivially  on  $C$: $C$  is
central in $\widehat G$ and Theorem \ref{maintheorem} is proved.

\newpage
\section{When $M$ is abelian} \label{proof2}

When $M$  is abelian, the  proof of Theorem \ref{maintheorem}  is more
involved.  We will  use  some  results from  \cite{V2}.  Since $M$  is
abelian and $P$  is a maximal parabolic subgroup, this  means that the
semisimple  part  of  $L$  has  $\Q$-rank  zero,  and  hence  $\Q-rank
(G)=1$. We state them now. \\

[1]  For each  prime  $p$,  denote by  $G_p  \subset  \widehat G$  the
subgroup generated by  $U^{\pm}(\Q _p)$.  The group $C$  is central if
and  only if,  for every  pair $p,q$  of distinct  primes, the  groups
$G_p,G_q$ commute. In \cite{V2} this is proved only in the case $C$ is
a (compact)  profinite group  (since the main  application was  to the
centrality of the congruence subgroup  kernel), but the proof works in
general and does not use the compactness of $C$. \\

[2] There  exists a  morphism $\phi: H=SL_2  \ra G$  of $\Q$-algebraic
groups such  that $\phi (U_H^{\pm})\subset U^{\pm}$.  Here $U_H^{\pm}$
is the  group of upper  (resp lower) triangular unipotent  matrices in
$SL_2$.  Further,  the  conjugates $\{^s\phi  (U_H^+),s  \in  L(\Q)\}$
generate the group $U^+$. \\

[3] There exists an infinite subgroup $\D \subset M(\Z)$ such that for
{\it every}  triple $a,b,k$  of mutually  coprime integers,  the group
generated  by the  collection $\{M(azk+b):  z \in  \Z\}$ of  subgroups
contains  this  fixed  group  $\D$,   and  such  that  the  commutator
$[\D,\phi (SL_2)]$ contains $\phi (SL_2)$. \\

Assume the  above facts,  and write $H=SL_2$.   For each  integer $k$,
write $F_H(k), \G _{H,k}$ for the  intersections $F(k)\cap H$ , $\G _k
\cap  H$.   Denote   by  $C_H(k)$  the  double   coset  $F_H(k)\bs  \G
_{H,k}/F_H(k)$ and  fix an element  $x \in C_H(k)$. Then  $F_H(k)$ has
the  same closure  in  the  congruence completion  of  $H(\Q)$ as  $\G
_{H,k}$.  If  $x=\begin{pmatrix} a & b  \\ c & d\end{pmatrix}$  and $c
\neq 0$ (which we may assume after replacing $x$ by a left translation
by  a   suitable  element  of   $F(k)$)  we  may  write   $x=vp$  with
$v=\begin{pmatrix} 1 &  0 \\ \frac{c}{a} & 1  \end{pmatrix}$ and hence
the common  denominator of $v,p$  is the  integer $a$. For  a suitable
power $N$ which depends only on the embedding $\phi :H \ra G$, we have
that   for  $m   \in   M(a^N)$,  the   commutator  $[m,\phi   (v)]=m\phi 
(v)m^{-1}\phi (v)^{-1} \in V^-(k\Z)  \subset F(k)$, and the commutator
$[(\phi (p))^{-1},m]\in (MV)(k\Z) \in  F(k)$. Consider the conjugate $
m\phi (x) m^{-1}$. Writing $x=vp$ we get
\[m\phi (x) m^{-1}=m\phi (v) m^{-1} m\phi (p) m^{-1}= \]
\[= [m, \phi (v)]x [\phi (p)^{-1},m] \in F(k)\phi (x) F(k).\] That is,
$m\phi (x)m^{-1}=\phi  (x)$ as double  cosets.  Thus the  group $M(a)$
fixes the  image of the  element $x \in C_H(k)$  under $\phi $  in the
double coset $C_k=F(k)\bs  \G_k/F(k)$ where $a$ is the  top left entry
of the matrix $x$. \\

We may replace  $x$ by an element $y=  x \g $ with $\g  \in F(k)$.  We
choose $\g  =\begin{pmatrix}1 & 0  \\ kz  & 1 \end{pmatrix}$  for some
integer  $z$.  Then  $y=\begin{pmatrix}  a+bzk  &  b  \\  c  +  dkz  &
d\end{pmatrix} $ has top left entry  $a+bkz$. By the conclusion of the
preceding paragraph, the  group $M(a+bkz)$ also fixes  the element $y$
viewed as  a double coset in  $C_k$.  But by construction  $x=y$ as
double cosets, and hence both the groups $M(a)$ and $M(a+bkz)$ fix $x$
for  all $z  \in  Z$.  Thus the  group  $M_{a,b,k}$  generated by  the
collection $\{ M(a+bkz)\Z): z\in \Z\}$ fixes the double coset $x$.  By
[3]  of the  listed  facts, there  is a  fixed  infinite subgroup  $\D
\subset M_{a,b,k}$ for  {\it every} $a,b,k$. Hence $\D$  fixes $x$ for
every $x\in  C_H(k)$ and hence $C_H(k)$  is fixed by $\D$  for every $k$.
By taking inverse limits, we see that the image of $C_H$ under the map
$\phi $ is fixed by all of $\D$.\\

However, the image $\phi (C_H)$ of $C_H$ is invariant under the action
of $SL_2=H$. Again by the second  part of [3], $\phi H\subset [\D,\phi
(H)]$  acts trivially  on $\phi  (C_H)$ and  hence $\phi  (C_H)$ is  a
central extension of $\phi (SL_2(\A_f))$.  Therefore, by fact [1], for
each pair of distinct primes $p,q$ the groups $\phi (U_H^+(\Q_p))$ and
$\phi (U^-(\Q_q))$ commute. \\

Let $s\in L(\Q)$  be arbitrary, and write  $s=(s_p)\in L(\A_f)$. Being
the linear action, the adjoint action of $L(\Q)$ on $U^{\pm}(\A_f)$ in
the topological group  $\widehat G$ factors through  the finite adelic
group $L(\A_f)$. Hence, for each $p,q$, $u \in U_H^+(\Q_p)$ and $v \in
U_H^-(\Q _q)$, we have
\[   s\phi   (u)   s^{-1}=s_p   \phi   (u)   s_p^{-1},   \quad   s\phi
(v)s^{-1}=s_q\phi  (v)s_q^{-1}.\] Furthermore,  by weak  approximation
(\cite{Gille}), $L(\Q)$ is  dense in $L(\Q _p)\times  L(\Q _q)$. Since
$\phi (u)$ and  $\phi (v)$ commute by the conclusion  of the preceding
paragraph,  we see  ( by  taking limits  of elements  in $L(\Q)\subset
L(\Q_p)\times L(\Q  _q)$) that for  every $s_p \in L(\Q_p)$  and every
$s_q \in L(\Q  _q)$, the elements $s_p\phi (u)  s_p^{-1}$ and $s_q\phi
(v) s_q^{-1}$ commute for  all $u,v$. But by fact [2],  $\phi $ may be
so  chosen  that the  collection  $s_p\phi  (u)s_p^{-1}$ with  $s_p\in
L(\Q_p),u\in U_H^+(\Q_p)$ generates all  of $U^+(\Q_p)$; similarly for
$U^-(\Q  _q)$. Hence  $U^+(\Q_p)$ commutes  with $U^-(\Q_q)$  for each
pair $p,q$  of distinct primes.  By  fact [1], this means  that $C$ is
central.  Thus we have proved Theorem \ref{maintheorem} in all cases. \\

Since  we  have  already  shown   in  the  introduction  that  Theorem
\ref{maintheorem}   implies   Theorems   \ref{commutatortheorem}   and
\ref{elementary}, we have also  proved Theorem \ref{elementary} in all
cases.

\end{document}